\documentclass[11pt]{article}

\usepackage{amsfonts}
\usepackage{amssymb}
\usepackage{eucal}
\usepackage[all]{xy}
\usepackage[small]{titlesec}

\date{}

{} 

\newtheorem{theorem}{{\bf Theorem}}[section]

\newtheorem{corollary}[theorem]{{\bf Corollary}}
\newtheorem{remark}[theorem]{{\bf Remark}} 
\newtheorem{problem}[theorem]{{\bf Problem}}

\newtheorem{proposition}[theorem]{{\bf Proposition}}

\newtheorem{definition}[theorem]{{\bf Definition}}
\newcommand{\proof}{{\bf Proof.  }}
\newcommand{\eproof}{{\square}}

\newcommand{\bN}{{\mathbb N}}

\title{\Large \bf A note on the Alster, Menger and $D$-type properties }

\author{Alexander V. Osipov\\
\vspace{-2mm}\footnotesize {Krasovskii Institute of Mathematics
and Mechanics, Ural Federal
 University, }\\
\footnotesize{Ural State University of Economics, Ekaterinburg, Russia}\\
\footnotesize{oab@list.ru}\\\\
Selma \"{O}z\c{c}a\u{g}\footnote{Corresponding author }\\
\vspace{-2mm}\footnotesize{ Department of Mathematics, Hacettepe
University,
Ankara, Turkey}\\
\footnotesize{sozcag@hacettepe.edu.tr}}

\newcommand{\sF}{\mathcal{F}}

\newcommand{\N}{{\mathbf N}}
\newcommand{\sone}{{\sf S}_1}
\newcommand{\gone}{{\sf G}_1}
\newcommand{\sfin}{{\sf S}_{fin}}
\newcommand{\gfin}{{\sf G}_{fin}}

\newcommand{\naturals}{{\mathbb N}}
\begin{document}

\maketitle

\begin{abstract}
In this paper we give  new characterizations for almost Menger and
weakly Menger spaces by neighborhood assignments  and  define a natural weakening of almost
$D$-spaces and weakly $D$-spaces. We  discuss the relationships  among the 
properties  ``D-type'', ``Menger'', ``Alster'', and the weak versions  of these properties in Lindel\"{o}f spaces. 

\end{abstract}

\begin{flushleft}
\vspace{.3cm} \noindent  {\sf Keywords}: Menger spaces, almost
Menger spaces, weakly Menger spaces, $D$-spaces, Alster spaces,
almost Alster spaces, selection principles,
topological games.
\end{flushleft}

\section{Introduction}

A neighborhood assignment for a topological space $(X,\tau)$ is a function
$\N:X\rightarrow \tau$ such that $x\in \N(x)$ for any $x\in X$.
A space $X$ is called a $D$-space if for any
neighborhood assignment $\N$ for $X$ there exists a closed
discrete subspace  $D$ of $X$ such that $\N(D)=\{\N(x):x\in D\}$
covers $X$, i.e., $X=\bigcup_{x \in D} \N(x)=\bigcup\N(D)$ \cite{dp}.
The $D$-property is a kind of covering property and it is easy to
see compact spaces and also $\sigma$-compact spaces are $D$-spaces,
however  it is still an open question whether
Lindel\"{o}f spaces are $D$-spaces.

\smallskip
$D$-spaces became the subject of many researches since it was first
introduced by E.K. van Dowen and W.F. Pfeffer in \cite{dp} and  many results
were obtained on the relations with the other covering
properties, generalizations and applications of $D$-spaces \cite{arh,au1,au2,auzd,eis,gr}.
There are two survey papers on $D$-spaces
\cite{eis,gr} and it is worth to mention that in \cite{yaj}  Szeptycki and Yajima desribed the most recent developments
 in the study of relationships of the  $D$-space to other covering properties.

\smallskip
Recently Aurichi in \cite{au1} proved that  Menger $T_1$ spaces are
$D$-spaces by using the partial open neighborhood assignment game
(PONAG) which shows that topological games are useful tools to
analyze the class of $D$-spaces.

\smallskip
Very recently Peng and Shen in \cite{peng} gave a new characterization for Menger spaces by neighborhood assignment and thus Aurichi's result \cite[Corollary 2.7]{au1} is easily obtained using this characterization.

\smallskip
In another direction the Menger property and its generalizations, almost Menger and weakly 
Menger properties have been investigated by many mathematicians (see \cite{babpan,koc,kocev,kocin,pansera}).

\smallskip
Now in this note we first give a new characterization for almost
Menger spaces by neighborhood assignment. We continue by introducing a new notion 
called ``open-dense almost(weakly) $D$-space'' and show that every Lindel\"{o}f almost Menger
space is an open-dense almost $D$-space. Similar results will be presented for the class of
weakly Menger spaces in Section $3$ and at the end of this section the relationships  among  all these covering spaces are investigated.

\smallskip
To conclude this  introduction we give some basic definitions  and a short summary for
selection principles and topological games. All topological spaces considered in this paper are assumed to be $T_1$. 

\smallskip
Recall that a topological space $(X,\tau)$ is Menger \cite{hur,menger}(almost
Menger, weakly Menger\cite{koc,kocin,pansera}) if for every sequence  $(\mathcal U_n :
n\in\naturals)$  of open covers of $X$, there exists a sequence
$(\mathcal V_n : n\in\naturals)$ of finite families such that for
each $n$, $\mathcal V_n\subseteq \mathcal U_n$ and
$\bigcup\{\mathcal V_n: n\in\naturals\}=X$(resp. $X=\bigcup_{n\in
\naturals}\{\overline {V}: V\in \mathcal V_n\}, \bigcup\{\mathcal
V_n :n\in\naturals\}$ is dense in $X$).

\smallskip
We use the following notations for collections of coverings of a topological space $X$:

\begin{itemize}
\item $\mathcal O$ is the
collection of all open covers of a space $X$;

\item $\overline{\mathcal O}$ is the collection of families $\mathcal U$ of open subsets of
the space for which $\{\overline{U}: U\in\mathcal U\}$ covers the space $X$;

\item  $\mathcal D$ is the collection of families of open
sets with union dense in the space $X$;

\item We denote by $\bN$ the set of positive integers and $\omega=\bN \cup \{0\}$.

\end{itemize}

Let us recall \cite{sch} the two classical selection principles in topological spaces.
Let $\mathcal A$ and $\mathcal
B$ be the collection of subsets of an infinite set $X$.

\smallskip
$\sfin({\mathcal A},{\mathcal B})$ denotes the hypothesis: for
each sequence $(A_n:n\in\naturals)$ of elements of ${\mathcal A}$
there is a sequence $(B_n:n\in\naturals)$ of finite sets such that
for each $n$, $B_n\subset A_n$, and $\bigcup_{n\in\naturals}B_n\in
{\mathcal B}$.

\smallskip
$\sone({\mathcal A},{\mathcal B})$ denotes the 
hypothesis: for each sequence $(A_n:n\in\naturals)$ of elements of
${\mathcal A}$ there is a sequence $(b_n:n\in\naturals)$ such that
for each $n$, $b_n\in A_n$, and $\{b_n:n\in\naturals\}$ is an
element of ${\mathcal B}$.

\smallskip
$\sfin({\mathcal O},{\mathcal O})$ denotes the Menger property,
$\sfin({\mathcal O},\overline{\mathcal O})$ denotes the almost Menger
property, and $\sfin({\mathcal O},{\mathcal D})$ denotes the
weakly Menger property.

\smallskip
There are games associated to the above selection principles. \\
The game  $\gfin({\mathcal A},{\mathcal B})$ associated with
$\sfin({\mathcal A},{\mathcal B})$ is as follows: Two players ONE and
TWO play an inning. In the $n$-th round  ONE chooses a set
$A_n\in{\mathcal A}$, while TWO responds by choosing a finite set
$B_n \subset A_n$. A play $(A_1,B_1,\cdots,A_n,B_n,\cdots)$ is won
by TWO if and only if $\bigcup_ {n\in \naturals}B_n \in {\mathcal
B}$. We denote the games with length $\omega$ by
$\gfin^{\omega}({\mathcal O},{\overline{\mathcal O}})$. In inning
$n<\omega$ ONE first selects an $A_n\in{\mathcal A}$, while TWO
responds by choosing a finite set $B_n \subset A_n$. A play
$(A_0,B_0, A_1,B_1,\cdots,A_n,B_n,\cdots)$ is won by TWO if
$\bigcup_ {n\in \naturals}B_n \in {\mathcal B}$.

\smallskip
It is evident that if ONE does not have a winning strategy in the game $\gfin({\mathcal A},{\mathcal B})$ (resp. $\gone({\mathcal A},{\mathcal B})$) then the selection hypothesis $\sfin({\mathcal A},{\mathcal B})$ (resp. $\sone({\mathcal A},{\mathcal B})$) is true. The converse implication need not be always true.

\smallskip
Babinkostova, Pansera and Scheepers  obtained in \cite[Theorem 29]{bab} that for  Lindel\"{o}f spaces,
$X$ is almost Menger if and only if   ONE does not have a winning strategy
in $\gfin^{\omega}({\mathcal O},{\overline{\mathcal O}})$.

\section{Almost Menger and open-dense almost D-spaces}
In this section we define the weaker  versions of $D$-spaces and discuss the relationships among the D-type-properties. We also give a new characterization for almost Menger spaces by neigborhood assignment.

\smallskip
The notion of almost $D$-space was introduced by Kocev in \cite{kocev}.
A topological space $(X,\tau)$ is an almost $D$-space
if for every function $\N:X\rightarrow \tau$ such that $x\in
\N(x)$, there exists a closed discrete subspace  $D$ of $X$ such
that  $\bigcup_{x \in D} \overline{\N(x)}=X.$

\medskip
Now we define a natural weakening of almost $D$-spaces.

\begin{definition}\label{wD}\rm  A topological space $(X,\tau)$ is a \emph{weakly $D$-space}
if for every function $\N:X\rightarrow \tau$ such that $x\in
\N(x)$, there exists a closed discrete subspace  $D$ of $X$ such
that  $\overline{\bigcup_{x \in D} \N(x)}=X$.
\end{definition}

\smallskip
We will now consider a new property called ``open-dense almost 
$D$-space (open-dense weakly $D$-space)''  and discuss the relationships with the 
other $D$-type properties.

\begin{definition}\label{opendense}\rm A topological space $(X,\tau)$ is an
\emph{open-dense $D$-space (open-dense almost $D$-space, open-dense
weakly $D$-space)} if for every function $\N: X\rightarrow \tau$
such that $x\in \N(x)$, there exists a closed discrete subspace
$D$ of $\bigcup_{x \in D} \N(x)$ such that $\bigcup_{x \in D}
\N(x)=X$ (resp. $\bigcup_{x \in D} \overline{\N(x)}=X$,
$\overline{\bigcup_{x \in D} \N(x)}=X$.)
\end{definition}

\begin{remark}\label{R}\rm
The following diagram indicates the relations among the $D$-type covering properties.
Here the words open-dense will be abbreviated to o-d. We note that $D$-spaces and open dense $D$-spaces are equivalent and the other implications in the diagram can be obtained easily by Definition \ref{wD} and Definition \ref{opendense}. 
\end{remark}

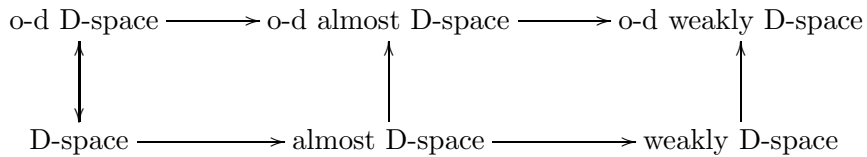
\[
\xymatrix @C=35pt @R=30pt{ \mbox{
o-d D-space}\ar[r]\ar[d] & \mbox{o-d almost D-space}\ar[r]&
\mbox{o-d weakly D-space}\\
\mbox{D-space}\ar[u]\ar[r]& \mbox{almost D-space}\ar[u]\ar[r] &
\mbox{weakly D-space} \ar[u]}
\]

\smallskip
In \cite{peng} the authors gave a new  characterization for Menger spaces 
by neighborhood assignments and thus they obtained Aurichi's  result ``Every 
$T_1$ Menger space is a $D$-space''  in connection with  this characterization. 

\smallskip
The following theorem gives a new  characterization for almost
Menger spaces by neighborhood assignment.

 \begin{theorem}\label{ana} Let $X$ be a Lindel\"{o}f space. Then $X$ is almost
 Menger if and only if for each sequence $\{\N_n : n\in\omega\}$ of
 neighborhood assignments for $X$, there exists for each
 $n\in\omega$ a finite subset $D_n$ of $X$ such that
 $\bigcup\{\overline{\N_n(x)} :x\in D_n, n\in \omega\}=X$ and
 $D=\bigcup\{D_n :n\in\omega\}$ is a closed discrete subspace of $\bigcup\{\N_n(x) :x\in D_n, n\in\omega\}$.
\end{theorem}

$\proof$ Let $X$ be an almost Menger space. The argument used here is
originally due to  Peng  and Shen by \cite[Theorem 3]{peng} and we will follow the similar lines
of that proof.

Let $(\N_n : n\in\omega)$ be a sequence of
 neighborhood assignments for $X$. The family $\{\N_0 (x):
 x\in X\}$ is an open cover of $X$ since $\N_0$ is a
 neighborhood assignment.

 \smallskip
 Let $\sigma$ be a strategy for player ONE and let $\mathcal
 U_0=\N_0=\sigma(\emptyset)$. Now let $\varphi$ be a way of player
 TWO. If $\varphi(\mathcal U_0)=\mathcal V_0$ then $\mathcal V_0
 \subset \mathcal U_0$ and $\mathcal V_0$ is a finite subfamily of
 $\mathcal U_0$. Then there is some finite subset $D_0 \subset X$
 such that $\mathcal V_0=\{\N_0 (x):
 x\in D_0\}$. Let $\mathcal U_1=\{\N_1(x)\cup\N_0(D_0) : x\in
 X\backslash \N_0(D_0)\}$ where $\N_0(D_0)=\bigcup\{\N_0(x) : x\in
 D_0\}$ then $\mathcal U_1$ is an open cover of $X$. Let
 $\mathcal U_1 =\sigma(\mathcal V_0)$. Then TWO responds with the
 finite subfamily $\mathcal V_1$ of $\mathcal U_1$. We denote
 $\mathcal V_1 =\varphi(\mathcal U_0,\mathcal U_1)$. So there is
 a finite set $D_1\subset X\backslash\N_0(D_0)$ such that
 $\mathcal V_1=\{\N_1 (x)\cup \N_0(D_0): x\in D_1\}$.

 \smallskip
 Now let $\mathcal U_2=\{\N_2(x)\cup \N_1(D_1)\cup\N_0(D_0) : x\in
 X\backslash (\N_0(D_0)\cup \N_1(D_1))\}$ and $\mathcal U_2$ is an
 open cover of $X$. Let
 $\mathcal U_2 =\sigma(\mathcal V_0,\mathcal V_1)$. This round  TWO responds with the
 finite subfamily $\mathcal V_2$ of $\mathcal U_2$ and denoted by $\mathcal V_2 =\varphi(\mathcal U_0,\mathcal U_1, \mathcal
 U_2)$. So there exists a finite set  $D_2\subset X\backslash\N_0(D_0)\cup \N_1(D_1)$ such that
 $\mathcal V_2=\{\N_2(x)\cup \N_1 (D_1)\cup \N_0(D_0): x\in D_2\}$.

\smallskip
The players proceed in this way and for $k\in\omega$ and
$k\geq 1$ we have\\  $(\mathcal U_0,\mathcal V_0; \mathcal
U_1,\mathcal V_1; \mathcal U_2,\mathcal V_2;...;\mathcal
U_k,\mathcal V_k)$ and a finite sequence $(D_i: i\leq k)$ of
finite subsets of $X$ satisfying the following conditions:

\begin{itemize}

\item $\mathcal U_i = \sigma(\mathcal V_0,...,\mathcal
V_{i-1})=\{\N_i(x)\cup((\cup\{\N_j(D_j): j<i\}):x\in X\backslash
\cup\{\N_j(D_j): j<i\}\}$ for each $0<i<k$ where
$\N_j(D_j)=\cup\{\N_j(x): x\in D_j\}$

\item $D_i\subset X\backslash \bigcup\{\N_j(D_j): j<i\}$ for each
$0\leq i \leq k$.

\item $\mathcal V_i=\{\N_i(x)\cup(\cup\{\N_j(D_j): j<i\}): x\in
D_i\}$ for each $0\leq i \leq k$.

\item $\mathcal V_i=\varphi(\mathcal U_0,...,\mathcal U_i)$ for each $0\leq i \leq k$.

\end{itemize}

If $\mathcal U_{k+1}=\{\N_{k+1}(x)\cup (\cup\{\N_j(D_j): j\leq
k\}): x\in X\backslash \cup\{\N_j(D_j): j\leq k\}\}$ then
$\mathcal U_{k+1}$ is an open cover of $X$.

\smallskip
Let $\mathcal U_{k+1} = \sigma(\mathcal V_0,...,\mathcal V_k)$.
Thus TWO responds with a finite subfamily $\mathcal V_{k+1}$ of
$\mathcal U_{k+1}$ and there is a finite subset $D_{k+1}\subset
X\backslash \cup\{\N_j(D_j): j\leq k\}$ such that $\mathcal
V_{k+1} = \{\N_{k+1}(x) \cup (\cup \{\N_j(D_j): j\leq k\}: x\in
D_{k+1}\}$.

\smallskip
If we proceed we have a play $(\mathcal U_0,\mathcal V_0; \mathcal
U_1,\mathcal V_1;...;\mathcal U_k,\mathcal V_k;...)$ such that
$\mathcal U_{k} = \sigma(\mathcal V_0,...,\mathcal V_{k-1})$ and
$\sigma$ is a strategy for ONE. Apply the fact that $X$ is almost
Menger ONE does not have a winning strategy in
$\gfin^{\omega}({\mathcal O},{\overline{\mathcal O}})$ by \cite[Theorem
29]{bab}.  So this strategy can not be a winning
strategy and $\varphi$ be a winning strategy for TWO. Thus  ONE
concludes that $\bigcup\{{\overline{\cup\mathcal V_k}} : k\in
\mathbb{N}\} = X$

\smallskip
Here $\mathcal V_k = \{\N_k(x)\bigcup (\cup \{\N_j(D_j): j< k\}):
x\in D_{k}\}$. $\bigcup \mathcal V_k = \bigcup \{\N_j(D_j) : j\leq
k\} = \bigcup \{\N_j(x): x\in D_j, j\leq k\}$ for each
$k\in\omega$ thus $\overline{\bigcup\mathcal V_k} = \bigcup
\{\overline{\N_j(x)} : x\in D_j, j\leq k\}$ for each
$k\in\omega$ .
 We have now $X=\bigcup \{\overline{\N_n(x)} : x\in D_n, n\in \mathbb{N}\}$ and $D_n$ is a finite subset of $X$ such
 that $D_n\subset X\backslash \bigcup\{\N_j(D_j): j< n\}$ for each
 $n\in\omega$. Hence, $D=\bigcup\{D_n :n\in\omega\}$ is a closed discrete subspace of $\bigcup\{\N_n(x) :x\in D_n, n\in\omega\}$.

\medskip
Now to prove the sufficiency: Let $(\mathcal U_n: n\in\omega)$
be a sequence of open covers of $X$ and let $n\in\omega$. For
each $x\in X$ there is some $U(n,x)\in\mathcal U_n$ such that
$x\in U(n,x)$. Let $\N_n(x)=U(n,x)$  for each $x\in X$. Thus $\N_n
= \{U(n,x): x\in X\}$ is a neighborhood assignment for $X$.

\smallskip
By the hypothesis there is a finite subset $D_n$ of $X$ such that
$X=\bigcup\{\overline{\N_n(x)} : x\in D_n, n\in \omega\}$,
$D=\bigcup\{D_n :n\in\omega\}$ is a closed discrete subspace of
$\bigcup\{\N_n(x) : x\in D_n, n\in \mathbb{N}\}$. Next set
$\mathcal V_n = \{\N_n(x):x\in D_n\}$ for each $n$. The sequence
$\mathcal V_n$ is a finite subfamily of $\mathcal U_n$ and $X =
\bigcup\{\overline{\bigcup\mathcal V_n} : n\in\omega\}$ so
that $X$ is almost Menger. $\eproof$

\begin{corollary}\label{cr1}
Every Lindel\"{o}f almost Menger space is an open-dense
almost $D$-space.
\end{corollary}

$\proof$ Let $\N$ be any neighborhood assignment for  $X$. For
each $n\in\omega$ we let $\N_n=\N$. So $(\N_n:n\in\omega)$
is a sequence of neighborhood assignments of $X$. By the Theorem
\ref{ana}, there exists a finite subset $D_n\subset X$ for each
$n\in\omega$ such that $X=\bigcup\{\overline{\N_n(x)} : x\in
D_n, n\in \omega\}$ and $D=\bigcup\{D_n :n\in\omega\}$ is a
closed discrete subspace of $\bigcup\{\N_n(x) : x\in D_n, n\in
\omega\}$. So $X=\bigcup\{{\overline{\N(x)}} : x\in D\}$. Thus
$X$ is open-dense almost $D$-space. $\eproof$

\medskip
We conclude this section by considering the relationship between almost Alster and open dense almost 
$D$-spaces.  Alster spaces were introduced to characterize the class of productively Lindel\"{o}f  spaces. We recall from \cite{alster} that a family $\sF$ of $G_\delta$ subsets of a space $X$ is a $G_\delta$ compact cover if 
for each compact subset $K$ of $X$ a set $F\in \sF$ such that $K\subset F$. $X$ is an Alster space if each $G_\delta$  compact cover of the space has a countable subset that covers $X$. We note that Alster spaces are Menger.

\smallskip
In \cite{kocev} Kocev  defined the  almost Alster spaces: A topological
space $X$ is almost Alster if each $G_\delta$-compact cover of $X$
has the countable subset $\mathcal V$ such that
$\bigcup\{{\overline {V}}:V\in\mathcal V\}=X$. It is shown that
almost Alster spaces are almost Menger \cite[Theorem 3.10]{kocev}.
Thus we immediately have the following result:

\begin{corollary}\label{cr2}
Let $X$ be a Lindel\"{o}f  space. Every almost Alster  space
is open-dense almost $D$-space.
\end{corollary}

\begin{remark} \rm The hypothesis that $X$ is Lindel\"{o}f in
Theorem \ref{ana}, Corollary \ref{cr1} and \ref{cr2} is not
necessary. Similarly to the proof of Theorem 12 in \cite{bab}, we
can prove that in some models of set theory there is a
non-Lindel\"{o}f space $X$ for which TWO has a winning
strategy in $G^{\omega}_{1}(\mathcal{O}, \overline{\mathcal{O}})$,
and, hence, in $G^{\omega}_{fin}(\mathcal{O},
\overline{\mathcal{O}})$.
\end{remark}

\begin{remark}\rm In \cite{kocev}  Kocev formulated a theorem \cite[Theorem~4.3]{kocev}:
 ``\emph{Let $X$ be a Lindel\"{o}f space. If $X$ is almost
Menger, then $X$ is an almost $D$-space.''} But Kocev proved only
that if $X$ is Lindel\"{o}f and almost Menger, then $X$ is an
open-dense almost $D$-space. We do not know whether Kocev's
hypothesis \cite[Theorem~4.3]{kocev} is correct.
\end{remark}

\section{Weakly Menger and open-dense weakly D-spaces}

In this section we  give the  characterization  of  weakly Menger spaces  as similar to characterization of almost Menger spaces as being done in previous section. 

\begin{theorem}\label{th2}
Let $X$ be a  Lindel\"{o}f  space. $X$ is weakly Menger, if
and only if for each sequence $(\N_n : n\in\omega)$ of
 neighborhood assignments for $X$ there exists for each
$n\in\omega$ a finite subset $D_n$ of $X$ such that
$\bigcup\{\N_n(x):x\in D_n, n\in \omega\}$ is a dense subset
of $X$ and  $D=\bigcup \{D_n: n\in\omega\}$ is closed discrete
subspace of $\bigcup\{\N_n(x):x\in D_n, n\in \omega\}$.
\end{theorem}

$\proof$ The proof proceed like  Theorem \ref{ana}
($\Rightarrow$). At one point in the proof of Theorem \ref{ana} we
use the game $\gfin^{\omega}({\mathcal O},{\overline{\mathcal
O}})$. Instead of using this game we apply the game
$\gfin^{\omega}({\mathcal O},\mathcal D)$ which means that there
is a winning strategy $\varphi$ for player TWO then
${\overline{\bigcup\{\bigcup\mathcal F_{k}: k\in \omega\}}}=
X$ and we obtain $X = {\overline{\bigcup\{N_n(x):x\in D_n, n\in
\omega\}}}$ for a closed discrete subspace $D$ of
$\bigcup\{\N_n(x):x\in D_n, n\in \omega\}$ where $D=\bigcup
\{D_n: n\in\omega\}$.

($\Leftarrow$) Let $(\mathcal U_n: n\in\omega)$ be a sequence
of open covers of $X$. $\N_n = \{\N_n(x):x\in X\}$ is  a
neighborhood assignment for $X$ where $\N_n(x)=U(n,x)$ and
$U(n,x)\in \mathcal U_n$ for each $x\in X$. Thus there is a finite
$D_n$ of $X$ such that $X = {\overline{\bigcup\{\N_n(x):x\in D_n,
n\in \omega\}}}$ and $D=\bigcup\{D_n: n\in\omega\}$ is a
closed discrete subspace of $\bigcup\{\N_n(x):x\in D_n, n\in
\omega\}$. Let $\mathcal V_n=\{\N_n(x):x\in D_n\}$ for each
$n$, then $\mathcal V_n$ is a finite subfamily of $\mathcal U_n$
and $X={\overline{\bigcup\{\bigcup\mathcal V_{n}: n\in
\omega\}}}$. $\eproof$

\medskip
The following result is clear from definitions and we omit the proof.

\begin{corollary}\label{cr3}
Every Lindel\"{o}f  weakly Menger  space is an open-dense
weakly $D$-space.
\end{corollary}

Now we will look at the relationship between weakly Alster spaces and 
open-dense weakly D-spaces. First we recall the notion of weakly Alster space which was defined in \cite{babpan}. Before giving the definition we need the following classes of covers of a topological space  \cite{babpan}: 

\smallskip
$\mathcal G_K$: The family consisting of sets $\mathcal U$ where
$X$ is not in $\mathcal U$, each element of $\mathcal U$ is a $
G_\delta$-set  and for each compact set $C\subset X$ there is
$U\in \mathcal U$ such that $C\subseteq U$.

\smallskip
$\mathcal G_D$: The collection of sets $\mathcal U$ where each
element of $\mathcal U$ is a $G_\delta$-set and $\bigcup \mathcal U$
is dense in the space.

\begin{definition} \rm\cite{babpan} A space $X$ is weakly Alster if each member of $\mathcal G_K$ has a countable subset
which is a member of $\mathcal G_D$.
\end{definition}

Note that a space $X$ is weakly Alster if and only if $X$ has the
property $S_1(\mathcal{G}_K, \mathcal{G}_D)$ \cite[Lemma 42]{babpan}.

\smallskip
The followings are  immediate and we omit the proofs.

\begin{proposition}
Every weakly Alster space is weakly Menger.
\end{proposition}

\begin{corollary}\label{cr4}
Let $X$ be a Lindel\"{o}f space. Every weakly Alster  space
is open-dense weakly $D$-space.
\end{corollary}

\begin{remark}\rm The hypothesis that $X$ is Lindel\"{o}f in
Theorem \ref{th2}, Corollary \ref{cr3} or \ref{cr4} is not
necessary: By Theorem 12 in \cite{bab} it is consistent that there
is a non-Lindel\"{o}f space $X$ for which TWO has a winning
strategy in $G^{\omega}_{1}(\mathcal{O}, \mathcal{D})$, and,
hence, in $G^{\omega}_{fin}(\mathcal{O}, \mathcal{D})$.
\end{remark}

We end this paper by  summarizing  the relationships between considered notions  for
Lindel\"{o}f spaces in the next diagram.

\bigskip
%\SelectTips{eu}{12}
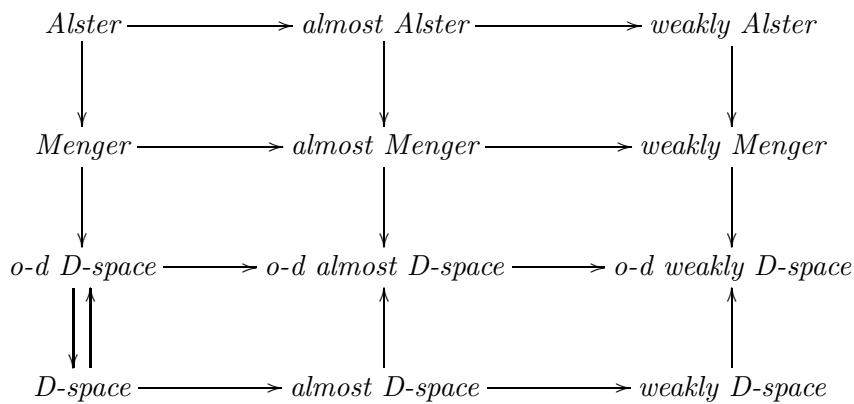
\[
\xymatrix @C=35pt @R=30pt{ \mbox{\it Alster}\ar[r]^-{}\ar[d]^{} &
\mbox{\it almost Alster}\ar[r]^-{}\ar[d]^{}&
\mbox{\it weakly Alster}\ar[d]^{}\\
\mbox{\it Menger}\ar[d]^{}\ar[r]^-{}& \mbox{\it almost
Menger}\ar[d]^{}\ar[r]^-{} &
\mbox{\it weakly Menger} \ar[d]^{}\\
\mbox{\it o-d D-space}\ar@<-.7ex>[d]\ar[r]^-{}& \mbox{\it o-d
almost D-space}\ar[r]^-{} &
\mbox{\it o-d weakly D-space}\\
\mbox{\it D-space}\ar@<-.7ex>[u]\ar[r]^-{}& \mbox{\it almost
D-space}\ar[u]_{}\ar[r]^-{ } & \mbox{\it weakly D-space}\ar[u]_{}
}
\]

\bigskip
\begin{remark}\rm Note that every Alster space is Hurewicz
\cite{tall}, and a Hurewicz space is Menger. A Menger subspace of
the real line which is not Hurewicz (is not Alster)
\cite{ts,tszd}. Note also that the space $\mathbb{P}$ of irrationals is Lindel\"{o}f and $D$-space but 
not Menger. In \cite{koc} it was shown that in regular spaces Menger property and almost Menger 
property are equivalent thus the space $\mathbb{P}$ of irrationals is a $D$-space but not almost Menger. 
\end {remark}

In \cite{bab} the authors asked whether there is a Lindel\"{o}f space which is not weakly Menger.  Sakai 
answered this question in the affirmative \cite [Theorem 2.5]{sakai}. 

\smallskip
\begin{problem} Is Sakai's example a $D$-space? (or may be a type of 
$D$-space).
\end{problem}

\end{document}